\documentclass{conm-p-l}

\usepackage[all]{xy}
\usepackage{graphics}

\theoremstyle{plain}
\newtheorem{theorem}{Theorem}[section]
\newtheorem{proposition}[theorem]{Proposition}
\newtheorem{lemma}[theorem]{Lemma}
\newtheorem{corollary}[theorem]{Corollary}
\newtheorem{conjecture}[theorem]{Conjecture}

\theoremstyle{definition}

\newtheorem{example}[theorem]{Example}
\newtheorem{question}[theorem]{Question}

\theoremstyle{remark}
\newtheorem{remark}[theorem]{Remark}

\newcommand{\thmref}[1]{Theorem~\ref{#1}}

\newcommand{\lemref}[1]{Lemma~\ref{#1}}
\newcommand{\corref}[1]{Corollary~\ref{#1}}
\newcommand{\conjref}[1]{Conjecture~\ref{#1}}
\newcommand{\remref}[1]{Remark~\ref{#1}}

\begin{document}

\title{The Rational Toomer invariant and Certain Elliptic Spaces}

\author{Gregory Lupton}
\address{Department of Mathematics,
         Cleveland State University,
         Cleveland OH 44115}
\email{Lupton@math.csuohio.edu}



\copyrightinfo{2002}{American Mathematical Society}

\begin{abstract}
We give an explicit formula for the rational category of an
elliptic space whose minimal model has a homogeneous-length
differential.  We also show that for such a space, there are no
gaps in the sequence of integers realized as the rational Toomer
invariant of some cohomology class.  With an additional
hypothesis, we show a result from which we deduce the relation
$\mathrm{dim}\big(H^*(X;\mathbb{Q})\big) \geq
2\,\mathrm{cat}_0(X)$.
\end{abstract}

\maketitle

\section{Introduction}

Since the fundamental paper of F{\'e}lix and Halperin
\cite{Fe-Ha82}, the so-called \emph{rational Toomer invariant} of
a space has played a central role in the development of results
concerning rational category. The rational Toomer invariant of a
space $X$ is a numerical rational homotopy invariant, denoted by
$\mathrm{e}_0(X)$, that provides a lower bound for the rational
category of a space. We recall its definition below.  In general,
its value is (strictly) between the rational cup length and  the
rational category of the space. One can also consider the rational
Toomer invariant of an individual cohomology class. This is a
finer invariant, whose supremum, taken over all cohomology classes
of a space, retrieves the rational Toomer invariant of the space.
The rational
 Toomer invariant  of a cohomology class can
also be identified with the strict, or essential, category weight
of the cohomology class, in the sense of Rudyak and Strom
(cf.~\cite{Rud99, Str97}). Since we are only concerned with the
rational case, we will henceforth refer to the rational Toomer
invariant simply as the Toomer invariant.

In this paper, we study the (rational) Toomer invariant  for
\emph{rationally elliptic} spaces.  Recall that a simply connected
space $X$ is called rationally elliptic if both
$H^*(X;\mathbb{Q})$ and $\pi_*(X)\otimes \mathbb{Q}$ are
finite-dimensional vector spaces. For these spaces, the Toomer
invariant  becomes all the more interesting, since it is known
that $\mathrm{cat}_0(X) = \mathrm{e}_0(X)$ for any elliptic space
$X$, where $\mathrm{cat}_0(X)$ denotes the \emph{rational
category} of the space $X$. Indeed, it is proved more generally in
\cite{F-H-L98} that $\mathrm{cat}_0(X) = \mathrm{e}_0(X)$ for any
space $X$ that satisfies rational Poincar{\'e} duality.

We now outline our results and indicate the organization of the
paper.  The main result is \thmref{thm:homogeneous, no e_0 gaps}.
We assume that the minimal model of an elliptic space has a
homogeneous-length differential. Then part (A) of
\thmref{thm:homogeneous, no e_0 gaps} gives a formula for the
Toomer invariant, and hence the rational category, of the elliptic
space.  The approach used to prove this part of the result leads
to finer information about the Toomer invariant  of individual
cohomology classes. Paraphrasing the precise statement, part (B)
of \thmref{thm:homogeneous, no e_0 gaps} says that for an elliptic
space whose minimal model has a homogeneous-length differential,
there are no gaps in the sequence of integers realized as the
Toomer invariant  of some cohomology class. Part (C) of the result
gives information about the location of cohomology classes that
have a given Toomer invariant. Our second result is
\thmref{thm:odd spherical generator}. Here we show that under a
fairly mild additional hypothesis, there are at least two linearly
independent cohomology classes with Toomer invariant  any integer
strictly between zero and the Toomer invariant  of the space.  The
inequality $\mathrm{dim}\big(H^*(X;\mathbb{Q})\big) \geq
2\,\mathrm{cat}_0(X)$ results.   The paper ends with a brief
section of examples and comments.

\smallskip

\textsc{A Little Archaeology and Acknowledgements:} The formula of
part (A) of \thmref{thm:homogeneous, no e_0 gaps} can be adapted
into a lower bound for the rational category of any elliptic
space---cf.~\remref{rem: inequality}. Indeed, an earlier version
of this paper included this as a separate result. During the
refereeing process, however, it was brought to my attention that
this lower bound appears as Corollary 3 in \cite{Gh-Je01}.
Further, the methods used there to obtain the lower bound, amongst
a number of other interesting related results, are comparable to
the methods of this paper.  A forthcoming paper \cite{Cu-Je01}
also contains a similar, but more general result. The referee
further pointed out to me that part (A) of
\thmref{thm:homogeneous, no e_0 gaps} appears as a special case of
Theorem 1 (Theorem 7) of \cite{Le-Mu01}. Their method of proof in
that paper, however, is substantially different from the one given
here and does not yield the finer information given by parts (B)
and (C) of \thmref{thm:homogeneous, no e_0 gaps}.  I thank for the
referee for bringing these articles to my attention.

\smallskip

We finish this introductory section with a brief review of some
ideas from rational homotopy theory.  All results of this paper
are proved using standard tools of the subject. We refer to
\cite{F-H-T01} for a general introduction to these techniques.  We
recall some of the notation here. By a \emph{minimal algebra} we
mean a free graded commutative algebra $\Lambda V$, for some
finite-type graded vector space $V$, together with a differential
$d$ of degree $+1$ that is \emph{decomposable}, i.e., satisfies $d
\colon V \to \Lambda^{\geq2}V$.  We assume that the minimal
algebra is simply connected, i.e., that the vector space $V$ has
no generators in degrees lower than $2$.  This assumption is
necessary in order to translate our algebraic results into
topological ones, although it is not strictly necessary for the
algebraic results themselves. A minimal algebra $\Lambda V$ is
called \emph{elliptic} if both $V$ and the cohomology algebra
$H^*(\Lambda V)$ are finite-dimensional vector spaces. If $\{ v_1,
\dots, v_n\}$ is a graded basis for $V$, then we write $\Lambda V$
as $\Lambda(v_1, \dots, v_n)$. A basis can always be chosen so
that $dv_1 = 0$ and for $i\geq 2$, $dv_i \in \Lambda(v_1, \ldots,
v_{i-1})$.

Every simply connected space with rational cohomology of
finite-type has a corresponding  \emph{minimal model}, which is a
minimal algebra that encodes the rational homotopy of the space.
Although our results are stated and proved in purely algebraic
terms, they do admit topological interpretations via this
correspondence.  In particular, a simply connected space is
rationally elliptic if its minimal model is an elliptic minimal
algebra.  Because of the correspondence between spaces and their
minimal models, this characterization of a rationally elliptic
space coincides with that given earlier.

We recall the definition, in minimal algebra terms, of the Toomer
invariant for a cohomology class and for a space. Suppose that
$\Lambda V$ is a minimal algebra. For $n\geq 1$, let $p_n$ denote
projection onto the quotient differential graded (DG) algebra
obtained by factoring out the DG ideal generated by monomials of
length at least $n+1$, thus
$$p_n \colon \Lambda V \to \frac{\Lambda V}{\Lambda^{\geq n+1} V}.$$
Define $\mathrm{e}_0(\Lambda V)$ to be the smallest $n$ such that
$p_n$ induces an injection on cohomology,  or set
$\mathrm{e}_0(\Lambda V) = \infty$ if there is no such smallest
$n$ (cf.~\cite[p.381]{F-H-T01}).  This can be extended to apply to
individual cohomology classes as follows: Suppose $x \in
H^*(\Lambda V)$ is some fixed, non-zero cohomology class. For
obvious degree reasons, $p_n^*(x) \not= 0$ for large enough $n$.
We define $\mathrm{e}_0(x)$ to be the smallest $n$ for which
$p_n^*(x) \not= 0$. If the set of integers $\{ \mathrm{e}_0(x)
\mid x \not=0 \in H^*(\Lambda V) \}$ has a maximum, then we see
that $\mathrm{e}_0(\Lambda V)$ is this maximum. Otherwise we have
$\mathrm{e}_0(\Lambda V) = \infty$.

Where convenient, we will use a standard observation about
$\mathrm{e}_0(\Lambda V)$ for a minimal algebra $\Lambda V$ whose
cohomology satisfies Poincar{\'e} duality. Namely, that if $\mu$
denotes a fundamental class of $H^*(\Lambda V)$, then
$\mathrm{e}_0(\Lambda V) = \mathrm{e}_0(\mu)$
\cite[Lem.10.1]{Fe-Ha82}.  It is well-known that an elliptic
minimal algebra satisfies Poincar{\'e} duality \cite[Th.3]{Hal77}.
By the \emph{formal dimension} of an elliptic minimal algebra
$\Lambda V$, we mean the largest $i$ for which $H^i(\Lambda V)
\not= 0$. If $x$ is an element of a graded vector space, then we
denote the degree of $x$ by $|x|$.

\section{Main Results}

We say that $\Lambda V$ has \emph{differential of
homogeneous-length} $l$ if $d \colon V \to \Lambda^{l} V$.  A
\emph{coformal} minimal algebra has differential of
homogeneous-length $2$.  In the homogeneous-length differential
case, the cohomology admits a second grading, $H^+(\Lambda V) =
\oplus_{k\geq1} H^*_k(\Lambda V)$, given by length of
representative cocycle.  Thus $x \not= 0 \in H^*_k(\Lambda V)$ if
and only if $\mathrm{e}_0(x) = k$. We extend this second grading
to include the degree zero component by setting $H^*_0(\Lambda V)$
equal to $\mathbb{Q}$ in degree zero, and zero elsewhere---see
(\ref{eqn:n_0=N_0}) below. Now introduce the following notation:
If $H^*_k(\Lambda V) \not= 0$, then set
$$n_k = \mathrm{min}\{i \mid H^i_k(\Lambda V) \not= 0\}
\quad \text{and} \quad N_k = \mathrm{max}\{i \mid H^i_k(\Lambda V)
\not= 0\}.$$
We make some observations about the bigraded Poincar{\'e} duality
algebra $H^*_*(\Lambda V)$.  Denote $\mathrm{e}_0(\Lambda V)$ by
$e$, so that $H^*(\Lambda V) = \oplus_{k = 0}^e H^*_k(\Lambda V)$.
By definition, we have
\begin{equation}\label{eqn:n_0=N_0}
H^i_0(\Lambda V) =
\begin{cases}
\mathbb{Q} & i = 0 \\
0 & i>0 \\
\end{cases}
\qquad \text{so} \quad n_0 = N_0 = 0.
\end{equation}
Now suppose that $H^*(\Lambda V)$  has formal dimension $N$. Since
the second grading comes from length of representative cocycle, we
have
\begin{equation}\label{eqn:n_e=N_e}
H^i_e(\Lambda V) =
\begin{cases}
\mathbb{Q} & i = N \\
0 & i < N \\
\end{cases}
\qquad \text{so} \quad n_e = N_e = N.
\end{equation}
We assemble our remaining remarks about $H^*_*(\Lambda V)$ into
the following lemma.

\begin{lemma}\label{lem:bigraded Pda} Let $H$ be a bigraded Poincar\'e duality
algebra of formal dimension $N$ that satisfies $H = \oplus_{k =
0}^e H^*_k$ together with (\ref{eqn:n_0=N_0}) and
(\ref{eqn:n_e=N_e}) above. Suppose that $H^*_k \not = 0$ for $k =
1, \ldots, e-1$, and that $p$ is some positive integer.  Then with
the above notation, the following are equivalent:
\begin{enumerate}

\item[(a)] $n_1 = p$ and $n_{k+1} \geq n_k + p$ for $k = 1, \ldots, e-1$.
\item[(b)] $N_{k+1} \geq N_{k} + p$ for $k = 0, \ldots, e-2$ and
$N = N_e = N_{e-1} + p$.
\end{enumerate}
\end{lemma}

\begin{proof} Because $H$ is a bigraded algebra, and also satisfies Poincar{\'e} duality,
it follows that there are non-degenerate pairings
$$H^{i}_{k}(\Lambda V) \times H^{N-i}_{e-k}(\Lambda V)
\to H^{N}_{e}(\Lambda V) \cong \mathbb{Q}\ ,$$
for $k = 1, \ldots, e-1$ and $i = 1, \ldots, N-1$. Hence $n_k = N
- N_{e-k}$ for $k = 1, \ldots, e-1$.  The equivalence of (a) and
(b) follows.
\end{proof}

To phrase the result, and to give its proof, we  say that a graded
vector space \emph{is $(i-1)$-connected} if degree $i$ is the
first non-zero degree. Also, for a bigraded cohomology algebra
$H^{*}_{*}(\Lambda V)$ we shall refer to elements in
$H^{i}_{*}(\Lambda V)$ as having \emph{upper degree} $i$ and to
elements in $H^{*}_{k}(\Lambda V)$ as having \emph{lower degree}
$k$.  In our setting, upper degree corresponds to the usual
topological degree, and lower degree corresponds to the Toomer
invariant  of a cohomology class. In this bigraded setting the
notation $|x|$ still denotes the upper degree of the element $x$.

\begin{theorem}\label{thm:homogeneous, no e_0 gaps} Suppose $\Lambda V$ is an elliptic
minimal algebra with homogeneous-length $l$ differential.  Set
$$e = \mathrm{dim}\,V^{\text{odd}} +
(l-2)\mathrm{dim}\,V^{\text{even}}.$$
Then we have
\begin{enumerate}
\item[(A)] $\mathrm{cat}_0(\Lambda V) = \mathrm{e}_0(\Lambda V) =
e$;
\item[(B)] $H^{*}_{k}(\Lambda V) \not= 0$ for each $k = 0, \ldots,
e$;
\item[(C)] Suppose that $V$ is $(p-1)$-connected.
 Then the bigraded Poincar{\'e} duality algebra
$H(\Lambda V)$ satisfies the two equivalent conditions of
\lemref{lem:bigraded Pda}.
\end{enumerate}
\end{theorem}

\begin{proof}
We will prove that $\mathrm{e}_0(\Lambda V) = e$. We can include
$\mathrm{cat}_0(\Lambda V)$ in the statement of part (A) due to
the result $\mathrm{cat}_0(\Lambda V) = \mathrm{e}_0(\Lambda V)$
for elliptic minimal algebras \cite{F-H-L98}.

Write $\Lambda V = \Lambda(x_1, \ldots, x_n)$, with $p=|x_1| \leq
|x_2| \leq \cdots \leq |x_n|$ and $d(x_1) = 0$.  We argue by
induction on the number of generators  $n$.  For $n=1$ ellipticity
requires that $|x_1|$ be odd. In this case $e=1$ and all parts of
the result are trivial.

Now assume inductively that the result holds for all elliptic
minimal algebras with homogeneous-length differential and fewer
than $n$ generators. Let $(\Lambda W; \bar d) = \Lambda(x_2,
\ldots, x_n; \bar d)$ be the quotient obtained by factoring out
the DG ideal generated by the generator $x_1$. Then $\Lambda W$ is
elliptic \cite[Prop.1]{Hal77} and has $n-1$ generators.   Note
that $W$ also has differential of homogeneous-length $l$. Both
$\Lambda V$ and $\Lambda W$ are bigraded DG algebras, with the
lower grading in each given by word-length.  Thus their cohomology
algebras are bigraded in the sense discussed above. We will make
extensive use of this bigraded structure without further remark.

Set $f = \text{dim}\,W^{\text{odd}} +
(l-2)\text{dim}\,W^{\text{even}}$, so that we have
\begin{equation}\label{eqn:f vs. e}
e =
\begin{cases}
f+1 & \text{if $|x_1|$ is odd} \\
f +  (l-2) &  \text{if $|x_1|$
is even.} \\
\end{cases}
\end{equation}
Let $M$ and $N$ denote the formal dimensions of $\Lambda W$ and
$\Lambda V$ respectively. Although we do not use it here, we
mention that there is a formula  that describes these formal
dimensions in terms of the degrees of the generators
\cite[Th.3]{Hal77}.

Further, set $m_k = \text{min}\{i \mid H^i_k(\Lambda W) \not= 0\}$
and $M_k = \text{max}\{i \mid H^i_k(\Lambda W) \not= 0\}$, for $k
= 0, \ldots, f$.   From the induction hypothesis, we have
$\mathrm{e}_0(\Lambda W) = f$, $H^*_k(\Lambda W) \not= 0$ for $k =
0, \ldots, f$ and $H^M_*(\Lambda W) = H^*_f(\Lambda W) =
H^M_f(\Lambda W) \cong \mathbb{Q}$. Also, $m_1 = |x_2|$ and
$m_{k+1} \geq m_k + |x_2|$ for $k = 1, \ldots, f - 1$.
Equivalently, as per \lemref{lem:bigraded Pda}, we have $M_{k+1}
\geq M_{k} + |x_2|$ for $k = 0, \ldots, f - 2$ and $M = M_f = M_{f
- 1} + |x_2|$. To prove the induction step there are two cases,
according as the parity of the degree of the first generator
$x_1$.

\smallskip

\textbf{Case I. $|x_1| = 2r+1$ is odd:}   In this case $e = f +
1$. Here we form the following short exact sequence of DG
\emph{vector spaces}:
$$\xymatrix{
0\ar[r]  &  \Lambda W \ar[r]^{j}& \Lambda V \ar[r]^{p}& \Lambda W
\ar[r] &0}.
$$
Here, $p$ denotes the projection and $j$ denotes the
degree-$(2r+1)$ linear map defined by
$$j(\chi) = (-1)^{|\chi|}\,x_1\,\chi\ $$
for $\chi \in \Lambda W$.   On cohomology, it is clear that $j^*$
increases lower degree by $1$, whilst $p^*$ preserves lower
degree. The corresponding long exact sequence in cohomology has
connecting homomorphism described as follows: Suppose $\chi \in
\Lambda W$ is any element. Regard $\chi$ as an element of $\Lambda
V$ and write $d\chi = \bar d\chi + x_1\theta(\chi)$. This defines
a  derivation $\theta \colon \Lambda W \to \Lambda W$, of
degree-$(-2r)$, that satisfies $\bar d \theta = \theta \bar d$.
Hence we have an induced derivation on cohomology $\theta^* \colon
H^{i}(\Lambda W) \to H^{i-2r}(\Lambda W)$, also of degree-$(-2r)$.
The derivation $\theta^*$ increases lower degree by $(l-2)$. The
resulting long exact cohomology sequence (the \emph{Wang
sequence}) is therefore a long exact sequence of bigraded
cohomology groups as follows:
\begin{equation}\label{eqn:long Wang bigraded}
\xymatrix{\cdots H^{i-1}_{k-1-(l-2)}(\Lambda W)
\ar[r]^-{\theta^*}& H^{i-2r-1}_{k-1}(\Lambda W) \ar[r]^-{j^*} &
H^{i}_{k}(\Lambda V) \ar[r]^-{p^*}& H^{i}_{k}(\Lambda W)\cdots}.
\end{equation}
It follows immediately from this long exact sequence that
$H^i_*(\Lambda V) = 0$ for $i \geq M+2r+2$ and  $H^*_k(\Lambda V)
= 0$ for $k \geq f+2 = e+1$. Furthermore, the long exact sequence
restricts to isomorphisms
\begin{equation}\label{eqn:Wang iso top dim bigraded}
\xymatrix{0 \ar[r]&  H^{M}_{k}(\Lambda W) \ar[r]^-{j^*}_-{\cong}&
H^{M+2r+1}_{k+1}(\Lambda V) \ar[r]&  0}.
\end{equation}
 One sees from this that the formal dimensions of
$\Lambda V$ and $\Lambda W$ are related by $N = M + 2r +1$ (cf.
\cite[Th.10.4]{Fe-Ha82}). This relation is also evident from the
formula of Halperin referred to above.  One also sees that
$\mathrm{e}_0(\Lambda V) = f+1 = e$. This proves part (A).

Now consider part (B).  In the previous paragraph, we showed that
$H^N_*(\Lambda V) = H^*_e(\Lambda V) = H^N_e(\Lambda V) \cong
\mathbb{Q}$.  It is automatic that $H^0_*(\Lambda V) =
H^*_0(\Lambda V) = H^0_0(\Lambda V) \cong \mathbb{Q}$. We consider
the remaining values $k = 1, \ldots, e-1$.  Denote the map
$$\theta^* \colon H^{*}_{k}(\Lambda W) \to H^{*}_{k+(l-2)}(\Lambda
W)$$
by $\theta^*_{k}$.  From the bigraded Wang sequence (\ref{eqn:long
Wang bigraded}) we obtain a short exact sequence
\begin{equation}\label{eqn:short Wang}
\xymatrix{0 \ar[r]&  \mathrm{cokernel}(\theta^*_{k-1-(l-2)})
\ar[r]^-{j^*}&
 H^{*}_{k}(\Lambda V) \ar[r]^-{p^*}& \mathrm{kernel}(\theta^*_{k}) \ar[r]&  0},
\end{equation}
for each $k = 1, \ldots, f$. From the inductive hypothesis applied
to $\Lambda W$, we have
\begin{equation}\label{eqn:Case I m}
m_{k+(l-2)} \geq m_{k} + (l-2) |x_2| \geq m_k
\end{equation}
and
\begin{equation}\label{eqn:Case II m}
M_{k+(l-2)} \geq M_{k} + (l-2) |x_2| \geq M_k.
\end{equation}
Since $\theta^*_k$ is of negative degree, (\ref{eqn:Case I m}) and
(\ref{eqn:Case II m}) imply respectively that both
$\text{kernel}(\theta^*_{k})$ and
$\text{cokernel}(\theta^*_{k-1-(l-2)})$ are non-zero for $k = 1,
\ldots, f$. It follows from (\ref{eqn:short Wang}) that
$\mathrm{dim}\big(H^{*}_{k}(\Lambda V)\big)\geq 2$ for each $k =
1, \ldots, f = e-1$.

Finally, we establish (C) in the present case, by showing that
$H^*_*(\Lambda V)$ satisfies condition (a) of \lemref{lem:bigraded
Pda} with $p = |x_1|$.  We observe from the sequence
(\ref{eqn:short Wang}) that $H^{i}_{k}(\Lambda V) = 0$ for $i <
m_{k-1} + |x_1|$. For if $i < m_{k-1} + |x_1|\leq m_{k-1} +
|x_2|$, it follows that there can be no contribution to
$H^{i}_{k}(\Lambda V)$ from
$\text{cokernel}\big(\theta^*_{k-1-(l-2)}\big)$. On the other
hand, $i < m_{k-1} + |x_1|$ implies $i < m_{k}$.   Thus if $i <
m_{k-1} + |x_1|$, then neither does $\text{kernel}(\theta^*_{k})$
contribute to $H^{i}_{k}(\Lambda V)$. So we have the inequality
\begin{equation}\label{eqn:nk lower}
n_{k} \geq m_{k-1} + |x_1|
\end{equation}
for $k = l, \ldots, f = e-1$. Next, recall that $\theta^*_k$ is of
negative degree. Hence $\text{kernel}(\theta^*_{k})$ begins in
degree $m_k$.  From (\ref{eqn:short Wang}), therefore, we obtain
\begin{equation}\label{eqn:nk upper}
n_{k} \leq m_{k}
\end{equation}
for $k = l, \ldots, f = e-1$. Combining the inequalities
(\ref{eqn:nk lower}) and (\ref{eqn:nk upper}), with a shift of
subscript in the latter, gives
$$n_{k} \geq m_{k-1} + |x_1| \geq n_{k-1} + |x_1|$$
for $k = 2, \ldots, f = e-1$.  In addition, it is evident that
$n_{1} = |x_1|$.  From (\ref{eqn:Wang iso top dim bigraded}), we
see that $n_{f+1} = m_{f} + |x_1| \geq n_{f} + |x_1|$, with the
latter inequality obtained from (\ref{eqn:nk upper}). This
establishes (C) and completes the induction step in Case I.

\smallskip

\textbf{Case II. $|x_1| = 2r$ is even:} In this case $e = f +
(l-2)$. Here, consider the short exact sequence (again of DG
vector spaces) as follows:
$$\xymatrix{
0\ar[r]  &  \Lambda V \ar[r]^{j}& \Lambda V \ar[r]^{p}& \Lambda W
\ar[r] &0},
$$
where $p$ denotes the projection and $j$ the degree-$(2r)$ map
defined by $j(\chi) = x_1\,\chi$ for $\chi \in \Lambda V$.  The
corresponding long exact sequence in cohomology has connecting
homomorphism as follows: Suppose $\bar d\chi = 0$ for $\chi \in
\Lambda W$, so that $d\chi = x_1\chi' = j(\chi')$ for some $\chi'
\in \Lambda V$ and  then $\partial^*([\chi]) = [\chi']$.  Here
$\partial^*$ is of degree-$(-2r+1)$ and if $\chi \in \Lambda^{\geq
k} W$, then $\chi' \in \Lambda^{\geq k+(l-2)} V$. On passing to
(bigraded) cohomology,  $\partial^*$ increases lower degree by
$(l-2)$, $p^*$ preserves lower degree and $j^*$ increases lower
degree by one.  In this case, therefore, the resulting long exact
sequence of bigraded cohomology groups (the \emph{Gysin sequence})
is as follows:
\begin{equation}\label{eqn:Gysin long bigraded}
\xymatrix{\cdots H^{i-2r}_{k-1}(\Lambda V) \ar[r]^-{j^*}&
H^{i}_{k}(\Lambda V)\ar[r]^-{p^*}& H^{i}_{k}(\Lambda W)
\ar[r]^-{\partial^*}& H^{i-2r+1}_{k+(l-2)}(\Lambda V)\cdots}.
\end{equation}
For $i \geq M+2$, this sequence gives isomorphisms
$$\xymatrix{0 \ar[r]& H^{i-2r}_*(\Lambda V)
\ar[r]^-{j^*}_-{\cong}&  H^{i}_*(\Lambda V)\ar[r]&  0}. $$
Since $H^{i}_*(\Lambda V)$ must be zero for sufficiently large
$i$, it follows that $N \leq M - 2r + 1$ and that $\partial^*
\colon H^{M}_{*}(\Lambda W) \to H^{M-2r+1}_{*}(\Lambda V)$ is an
isomorphism. Thus we obtain $N = M - 2r + 1$ (cf.
\cite[Th.10.4]{Fe-Ha82}). Once again, this relation is also
evident from the formula of Halperin referred to above. On the
other hand, for $k \geq f + (l-2) +2$, (\ref{eqn:Gysin long
bigraded}) gives isomorphisms
$$\xymatrix{0 \ar[r]& H^{*}_{k-1}(\Lambda V)
\ar[r]^-{j^*}_-{\cong}&  H^{*}_{k}(\Lambda V)\ar[r]&  0}. $$
Again, since $H^{*}_{k}(\Lambda V)$ must also be zero for
sufficiently large $k$, it follows that $\mathrm{e}_0(\Lambda V)
\leq f + (l-2) = e$. The isomorphism $\partial^*$ restricts to
isomorphisms
$$\xymatrix{0 \ar[r]^-{p^*}& H^{M}_{k}(\Lambda W)
\ar[r]^-{\partial^*}_-{\cong}& H^{M-2r+1}_{k+(l-2)}(\Lambda
V)\ar[r]^-{j^*}& 0}.
$$
It follows that $\mathrm{e}_0(\Lambda V) = f + (l-2) = e$, as
desired for part (A).

For parts (B) and (C), we show $H^{*}_{k}(\Lambda V) \not= 0$ for
$k = 1, \ldots, e$ and establish condition (b) of
\lemref{lem:bigraded Pda}, by using a (secondary) induction on
$k$.

The secondary induction hypothesis is as follows:
\begin{enumerate}
\item[(i)] $H^{*}_{i}(\Lambda V)\not=0$ for $i = 1, \dots, k-1$,
\item[(ii)] $N_i \geq N_{i-1} + 2r$ for $i = 1, \dots, k-1$, and
\item[(iii)] $N_i \geq M_{i - (l-2)} -2r +1$, for $i = l-1,
\dots, k-1$.
\end{enumerate}
This induction starts with $k=l$. Since the differential is of
length $l$, there are no boundaries of length $\leq (l-1)$.
Therefore $x_1^i \not=0 \in H^{2ri}_i(\Lambda V)$ for $i = 1,
\dots, l-1$.  Furthermore, if $\alpha \not=0 \in H^{N_i}_i(\Lambda
V)$ for $i = 1, \dots, l-2$, then $x_1\alpha \not=0 \in
H^{N_i+2r}_{i+1}(\Lambda V)$, again because there are no
boundaries of  length $\leq (l-1)$.  Hence $N_i \geq N_{i-1} + 2r$
for $i = 1, \dots, l-1$. This establishes (i) and (ii) with $k=l$.
For (iii), we must show that $N_{l-1} \geq M_{1} -2r +1$. So
consider the following portion of the Gysin sequence:
\begin{equation}\label{eqn:Gysin induction starts}
\xymatrix{H^{*}_{1}(\Lambda
V)\ar[r]_{p^*}\ar@/^2pc/[rr]^-{(j^*)^{l-2}}& H^{*}_{1}(\Lambda
W)\ar[r]_{\partial^*}& H^{*}_{l-1}(\Lambda V)\ar[r]_{j^*}&
H^{*}_{l}(\Lambda V)}.
\end{equation}
The map $(j^*)^{l-2} \colon H^{*}_{1}(\Lambda V) \to
H^{*}_{l-1}(\Lambda V)$ is simply $j^*$ composed with itself
$(l-2)$-times.  In other words, this map is multiplication by the
class of $x_1$ $(l-2)$-times.  Note that the diagram as displayed
is \emph{not} commutative.  By our primary induction hypothesis,
we have some $\alpha \not=0 \in H^{M_1}_1(\Lambda W)$. If
$\partial^*(\alpha) \not=0 \in H^{M_1-2r+1}_{l-1}(\Lambda V)$,
then $N_{l-1} \geq M_{1} -2r +1$, and (iii) holds for $k = l$.  On
the other hand, if $\partial^*(\alpha) =0$,  then exactness
provides some $\beta \not=0 \in H^{M_1}_1(\Lambda V)$ with
$p^*(\beta) = \alpha$.   Then $(j^*)^{l-2}(\beta) \not=0\in
H^{M_1+2r(l-2)}_{l-1}(\Lambda V)$, once more because there are no
boundaries of length $\leq l-1$. Since $|(j^*)^{l-2}(\beta)| =
M_1+2r(l-2) > M_1 - 2r +1$, we have shown (iii) for $k = l$ under
either of the two possible circumstances. This starts our
induction.

Now we make the inductive step.  Given the secondary inductive
hypothesis for some $k \geq l$, we will show (i), (ii) and (iii)
for $i = k$. Observe, that in any one Gysin sequence, some terms
cannot appear due to the non-consecutive nature of the indices. We
take the following two versions of the Gysin sequence, each of
which features $H^*_k(\Lambda V)$:
$$
\xymatrix{ & H^{*}_{k-1-(l-2)}(\Lambda V)\ar[dd]^{p^*} &  &
H^{*}_{k-1}(\Lambda V)\ar@/^1pc/[r]^-{j^*} &  H^{*}_{k}(\Lambda V)\ar[dd]^{p^*}  \\
 \ar@{.>}[ur] & & & & &\\
& H^{*}_{k-1-(l-2)}(\Lambda W)\ar[rruu]^{\partial^*} & & &
H^{*}_{k}(\Lambda W) \ar@{.>}[ur] }
$$
and
$$
\xymatrix{ & H^{*}_{k-(l-2)}(\Lambda V)\ar[dd]^{p^*} &  &
H^{*}_{k}(\Lambda V)\ar@/^1pc/[r]^-{j^*} &  H^{*}_{k+1}(\Lambda V)\ar[dd]^{p^*}  \\
 \ar@{.>}[ur] & & & & &\\
 & H^{*}_{k-(l-2)}(\Lambda W)\ar[rruu]^{\partial^*} &   & &  H^{*}_{k+1}(\Lambda W) \ar@{.>}[ur] }
$$
We splice these together at the $H^{*}_{k}(\Lambda V)$ term, then
 add the map $(j^*)^{l-2} = j^*\circ(j^*)^{l-3}\colon
H^{*}_{k-(l-2)}(\Lambda V) \to H^{*}_{k}(\Lambda V)$, as we did in
(\ref{eqn:Gysin induction starts}) when starting this induction.
This gives the following diagram, that contains the maps to which
we refer in our argument.
$$
\xymatrix{\cdots H^{*}_{k-1-(l-2)}(\Lambda V)\ar@/^2pc/[r]^-{j^*}
\ar[dd]_(0.3){p^*}  & H^{*}_{k-(l-2)}(\Lambda
V)\ar@/^2pc/[r]^-{(j^*)^{l-3}} \ar'[d]_{p^*}[dd] &
H^{*}_{k-1}(\Lambda V)\ar@/^2pc/[r]^-{j^*}\ar@{.>}'[d]_{p^*}[dd] &
  H^{*}_{k}(\Lambda V)\ar@{.>}[dd]_(0.3){p^*}\cdots\\
  & & \\
H^{*}_{k-1-(l-2)}(\Lambda W) \ar[rruu]^(0.3){\partial^*}&
H^{*}_{k-(l-2)}(\Lambda W)\ar[rruu]^(0.3){\partial^*} & &
 }
$$
We note again that this is \emph{not} a commutative diagram.

Our primary inductive hypothesis gives some $\alpha \not=0 \in
H^{M_{k-(l-2)}}_{k-(l-2)}(\Lambda W)$. If $\partial^*(\alpha)
\not=0 \in H^{M_{k-(l-2)}-2r+1}_{k}(\Lambda V)$, then $N_{k} \geq
M_{k-(l-2)} -2r +1$, and (iii) holds for $i=k$.  On the other
hand, if $\partial^*(\alpha) =0$,  then exactness provides some
$\beta \not=0 \in H^{M_{k-(l-2)}}_{k -(l-2)}(\Lambda V)$ with
$p^*(\beta) = \alpha$. We claim that $(j^*)^{l-2}(\beta) \not=0\in
H^{M_{k-(l-2)}+2r(l-2)}_{k}(\Lambda V)$. This claim follows from a
combination of exactness and degree considerations. For the case
$l=2$, we interpret $(j^*)^{0}(\beta)$ as $\beta$.  For $l\geq 3$,
Consider the sequence of elements $\{ \beta, j^*(\beta), \dots,
(j^*)^{t}(\beta), \dots, (j^*)^{l-2}(\beta)\}$. By exactness, we
have
$$
\begin{aligned}
\mathrm{ker}\big((j^*) \colon &H^{*}_{k-(l-2)+(t-1)}(\Lambda V)
\to H^{*}_{k-(l-2)+t}(\Lambda V)\big) \\
&= \mathrm{im}\big((\partial^*) \colon
H^{*}_{k-2(l-2)+(t-1)}(\Lambda W) \to
H^{*}_{k-(l-2)+(t-1)}(\Lambda V)\big)
\end{aligned}
$$
for $t=1, \dots, l-2$.    Now
$\partial^*(H^{*}_{k-2(l-2)+(t-1)}(\Lambda W))$ is zero in upper
degrees above degree $M_{k-2(l-2)+(t-1)} - 2r +1$.  Since
$|(j^*)^{t-1}(\beta)| = M_{k-(l-2)}+2r(t-1)$, and since our
primary induction hypothesis gives
$$
M_{k-(l-2)} \geq M_{k-2(l-2)+(t-1)} + \big((l-2) - (t-1)
\big)|x_2|,
$$
it follows that $(j^*)^{t-1}(\beta)$ is not in $\mathrm{ker}(j^*)$
for $t =1, \dots, l-2$. From the claim, we have  $N_k \geq
M_{k-(l-2)}+2r(l-2)
> M_{k-(l-2)} - 2r +1$. Under either of the two possible circumstances, therefore, we have
shown (i) and (iii) for $i=k$.

We complete the inductive step by showing (ii) for $i = k$. Our
secondary induction hypothesis provides an element $\gamma \not=0
\in H^{N_{k-1}}_{k-1}(\Lambda V)$. If $j^*(\gamma) \not=0 \in
H^{N_{k-1}+2r}_{k}(\Lambda V)$, then $N_{k} \geq N_{k-1} +2r$, and
(ii) holds for $i=k$.  On the other hand, if $j^*(\gamma) =0$,
then exactness implies $\gamma = \partial^*(\delta)$, for some
$\delta \not=0 \in H^{N_{k-1}+2r-1}_{k-1 -(l-2)}(\Lambda W)$. This
implies $M_{k-1 -(l-2)} \geq N_{k-1}+2r-1$ which, combined with
(iii) of the secondary induction hypothesis gives $M_{k-1 -(l-2)}
= N_{k-1}+2r-1$. Combining this equality with (iii) for $i=k$,
which we just established, and the primary induction hypothesis,
we have
$$
\begin{aligned}
N_k & \geq M_{k-(l-2)} -2r +1  \\
&\geq  M_{k-(l-2)} + |x_2| -2r +1 \\
& = N_{k - 1} + |x_2| \geq N_{k - 1} + |x_1|.\\
\end{aligned}
$$
This establishes (ii) in the case $i=k$ and hence completes our
secondary induction step.

To complete the primary induction step in Case II, it only remains
to show $N_e = N_{e-1} + |x_1|$. To see this, observe that
$H^*_*(\Lambda V)$ is a bigraded Poincar\'e duality algebra that
satisfies, in particular, $n_1 = N - N_{e-1}$. Clearly $n_1 =
|x_1|$ and so $N = N_{e-1} + |x_1|$. This completes the (primary)
induction step in Case II.

In both Case I and Case II the induction has been completed and
this proves the result.
\end{proof}

\begin{remark}
In the above, we referred to (\ref{eqn:long Wang bigraded}) as a
Wang sequence and (\ref{eqn:Gysin long bigraded}) as a Gysin
sequence. In Case I of the proof, we have a fibration sequence of
minimal models $\xymatrix{\Lambda(x_1)\ar[r]& \Lambda V \ar[r]&
\Lambda W}$. Topologically, this corresponds to a fibration
sequence of spaces $F \to E \to S^{2r+1}$, with base an
odd-dimensional sphere.  The Wang sequence of this fibration, in
the usual sense \cite[p.319]{Whi78}, corresponds to (\ref{eqn:long
Wang bigraded}).  In Case II, our fibration sequence of minimal
models corresponds topologically to a fibration sequence of spaces
$F \to E \to K(\mathbb{Q}, 2r)$.  Going one step back in the fibre
sequence gives a fibration sequence $K(\mathbb{Q}, 2r-1) \to F \to
E$. Since $K(\mathbb{Q}, 2r-1)$ and $S^{2r-1}$ have the same
rational homotopy type, we end up rationally with a fibration
sequence $S^{2r-1} \to F \to E$, with fibre an odd-dimensional
sphere.  The Gysin sequence of this fibration \cite[p.357]{Whi78},
corresponds to (\ref{eqn:Gysin long bigraded}).
\end{remark}

\begin{remark}\label{rem: inequality}
In the introduction, we mentioned that the formula in part (A) of
\thmref{thm:homogeneous, no e_0 gaps} can be adapted into a lower
bound for the rational category of any elliptic space.
Specifically, this goes as follows: We say $\Lambda V$ has
\emph{differential of length at least $l$} if $d \colon V \to
\Lambda^{\geq l} V$.  Then with $e$ as in the statement of
\thmref{thm:homogeneous, no e_0 gaps}, the inequality
$\mathrm{cat}_0(\Lambda V) = \mathrm{e}_0(\Lambda V) \geq e$ holds
for any elliptic space. Once again, this is Corollary 3 of
\cite{Gh-Je01}.  The proof of \thmref{thm:homogeneous, no e_0
gaps} is easily adapted to establish this inequality.  One simply
omits all reference to parts (B) and (C) of the proof, and uses
the Wang and Gysin sequences in their ordinary, i.e., not
bigraded, versions in precisely the same way as they were used
above.  Likewise, the proof of part (A) of
\thmref{thm:homogeneous, no e_0 gaps} can be given independently
of the proof of parts (B) and (C).  It is interesting to note,
however, that the proofs of these latter two parts cannot be
separated from each other.
\end{remark}

For our last result, we add the hypothesis that $\mathrm{ker}( d
\colon V^{\mathrm{odd}} \to \Lambda V)$ is non-zero to those of
\thmref{thm:homogeneous, no e_0 gaps}.  This extra hypothesis
admits a simple topological translation. Namely, that the rational
Hurewicz homomorphism $h \colon \pi_*(X)\otimes\mathbb{Q} \to
H_*(X;\mathbb{Q})$ is non-zero in some odd degree.

\begin{theorem}\label{thm:odd spherical generator} Suppose $\Lambda V$ is an elliptic
minimal algebra with homogeneous-length $l$ differential and
$\mathrm{ker}( d \colon V^{\mathrm{odd}} \to \Lambda V)$ non-zero.
Then $\mathrm{dim}\big(H^{*}_{k}(\Lambda V)\big)\geq 2$ for each
$k = 1, \ldots, e-1$, where $e = \mathrm{cat}_0(\Lambda V) =
\mathrm{e}_0(\Lambda V)$ is given by the formula of
\thmref{thm:homogeneous, no e_0 gaps}.
\end{theorem}

\begin{proof}
Suppose $u \in V^{\mathrm{odd}}$ is a cocycle.  Then we can write
$\Lambda V = \Lambda(x_1, \ldots, x_n)$ with $x_1 = u$. Now argue
exactly as in case I of the proof of \thmref{thm:homogeneous, no
e_0 gaps}. Observe that our quotient $\Lambda W =  \Lambda(x_2,
\ldots, x_n; \bar d)$ in the present result satisfies the
hypotheses of \thmref{thm:homogeneous, no e_0 gaps}.  Therefore,
the part of the argument featuring the short exact sequence
(\ref{eqn:short Wang}), with the inequalities (\ref{eqn:nk lower})
and (\ref{eqn:nk upper}), used there to establish the inductive
step, can be used here to conclude our result.
\end{proof}

\begin{corollary}\label{cor:dimH>=2cat}
Let $X$ be an elliptic space whose minimal model has a homo-
geneous-length differential and whose rational Hurewicz
homomorphism is non-zero in some odd degree. Then
$\mathrm{dim}\big(H^*(X;\mathbb{Q})\big) \geq 2\,\mathrm{cat}_0(X)
= 2\,\mathrm{e}_0(X)$.
\end{corollary}

\begin{proof}
Let $\Lambda V$ denote the minimal model of $X$.  From
\thmref{thm:odd spherical generator}, we obtain
$\mathrm{dim}\big(H^*(\Lambda V)\big) \geq
2\,\mathrm{cat}_0(\Lambda V) = 2\,\mathrm{e}_0(\Lambda V)$.  The
statements about $X$ follow from the standard translation from
minimal models to spaces.
\end{proof}

\section{Examples and Comments}

There are a number of special cases of our main results, which
either retrieve well-known results or provide interesting
examples.

\begin{example} Consider the case in which $l=2$.
Part (A) of \thmref{thm:homogeneous, no e_0 gaps} specializes to
retrieve part of \cite[Prop.10.6]{Fe-Ha82}, where it is shown that
$\mathrm{e}_0(\Lambda V) = \text{dim}\,V^{\text{odd}}$ for
$\Lambda V$ elliptic and coformal.
\end{example}

\begin{example}
Consider the case in which $\Lambda V = \Lambda(x_1, \ldots, x_n)$
is a minimal algebra with $|x_i| =1$ for each $i$. This type of
example arises as the minimal model of a \emph{nilmanifold}
\cite{Opr92}. Although we are primarily interested in the simply
connected case, the results proved here carry though verbatim for
this nilmanifold case.  Here $\mathrm{e}_0(x) = i$ if and only if
$x \not=0 \in H^i(\Lambda V)$.  Thus the Toomer invariant of a
cohomology class is identified with the degree and an integer $i$
is realized as the Toomer invariant  of some class if and only
$b_i \not=0$, where $b_i$ denotes the $i$'th Betti number of
$\Lambda V$, or of the nilmanifold of which $\Lambda V$ is the
minimal model.

For degree reasons, the differential here must be homogeneous of
length $2$.  Therefore \thmref{thm:odd spherical generator} and
\corref{cor:dimH>=2cat} specialize to yield the following
well-known result, which is essentially due to Dixmier
\cite{Dix55}:
\begin{proposition} A nilmanifold $X$ of dimension $n$ has $b_i
\geq 2$ for $1 \leq i \leq n-1$ and hence
$\mathrm{dim}\big(H^*(X;\mathbb{Q})\big) \geq 2
\,\mathrm{dim}(X)$.
\end{proposition}
\end{example}

The conclusion of \thmref{thm:odd spherical generator} obviously
holds in many cases besides those covered by the hypotheses. We
cannot resist making the following conjecture:

\begin{conjecture}\label{conj:dim Hi geq 2} Let $\Lambda V$ be elliptic
with homogeneous-length differential. Either
$\mathrm{dim}\big(H^{*}_{k}(\Lambda V)\big)\geq 2$ for each $k =
1, \ldots, e-1$, where $e = \mathrm{cat}_0(\Lambda V) =
\mathrm{e}_0(\Lambda V)$ is given by the formula of
\thmref{thm:homogeneous, no e_0 gaps}, or $H^{*}(\Lambda V)$ is a
truncated polynomial algebra on a single generator.
\end{conjecture}

In view of \thmref{thm:odd spherical generator}, one need only
consider the case in which $\mathrm{ker}( d \colon V \to \Lambda
V)$ is concentrated in even degrees.  We believe that an argument
as in case II of the proof of \thmref{thm:homogeneous, no e_0
gaps}, together with a careful analysis of the exceptional cases,
will establish \conjref{conj:dim Hi geq 2} at least in the
coformal case, $l = 2$.  We have been unable to prove the general
case, however, using this approach.

It is also clear that in many cases to which \thmref{thm:odd
spherical generator} applies, and others to which it does not, the
inequalities  $\mathrm{dim}\big(H^{*}_{k}(\Lambda V)\big)\geq 2$
for each $k = 1, \ldots, e-1$ and
$\mathrm{dim}\big(H^*(X;\mathbb{Q})\big) \geq
2\,\mathrm{cat}_0(X)$ are by no means sharp.  The following
examples show that, nonetheless, these inequalities are best
possible for a general result of this nature.

\begin{example}
Let $X_l = \mathbb{C}P^{l-1}\times S^{2r+1}$ for $l \geq 2$. Then
$X_l$ is elliptic and has minimal model with homogeneous-length
$l$ differential.  As is well known,  $\mathrm{cat}_0(X_l)=
\mathrm{e}_0(X_l) = l+1$. On the other hand, we see that
$\mathrm{dim}\big(H_k^*(X;\mathbb{Q})\big) = 2$ for $k = 1,
\ldots, l$.  Thus we have
$\mathrm{dim}\big(H^*(X_l;\mathbb{Q})\big) =
2\,\mathrm{cat}_0(X_l)$.
\end{example}

\begin{example}
Let $\Lambda V = \Lambda(x_1, x_2, y_1, y_2, y_3)$ with $|x_1| =
|x_2| = 2$, $|y_1| = |y_2| = |y_3| = 3$  and differential $d(x_1)
= 0 = d(x_2)$, $d(y_1) = x_1^2$, $d(y_2) = x_1x_2$ and $d(y_3) =
x_2^2$.   Then $\Lambda V$ is elliptic and coformal.  As is easily
checked, we have $\mathrm{cat}_0(\Lambda V)= \mathrm{e}_0(\Lambda
V) = 3$, and $\mathrm{dim}\big(H_k^*(\Lambda V)\big) = 2$ for $k =
1, 2$. In particular, we have $\mathrm{dim}\big(H^*(\Lambda
V)\big) = 2\,\mathrm{cat}_0(\Lambda V)$.
\end{example}

It would be interesting to know whether there are other examples
in which the inequalities of \corref{cor:dimH>=2cat} and
\thmref{thm:odd spherical generator} are sharp, or whether these
examples are essentially the only such.  In particular, it would
be interesting to find general conditions under which the
inequality of \corref{cor:dimH>=2cat} could be (substantially)
strengthened.

Finally, we remark that the original motivation for this work came
from a question of Yves F{\'e}lix, as to whether there exists any
space with ``$\mathrm{e}_0$-gaps'' in its cohomology. Precisely,
we say that a space $X$ has an ``$\mathrm{e}_0$-gap'' in its
cohomology if $H^*(X, \mathbb{Q})$ has an element whose Toomer
invariant  is $k$, but does not have any element whose Toomer
invariant  is $k-1$.   A recent example due to Kahl and
Vandembroucq \cite{KaVa01} shows that $\mathrm{e}_0$-gaps can
occur in the cohomology of a finite complex.  Their example is
actually a Poincar{\'e} duality space, but it is hyperbolic and
not elliptic. This leaves the following question:

\begin{question}
Can an \emph{elliptic} space have $\mathrm{e}_0$-gaps in its
cohomology?
\end{question}

If it is not possible for an elliptic space to have
$\mathrm{e}_0$-gaps in its cohomology, then it would seem
reasonable to extend  \conjref{conj:dim Hi geq 2} to the general
elliptic space.

\providecommand{\bysame}{\leavevmode\hbox
to3em{\hrulefill}\thinspace}
\providecommand{\MR}{\relax\ifhmode\unskip\space\fi MR }
\providecommand{\MRhref}[2]{%
  \href{http://www.ams.org/mathscinet-getitem?mr=#1}{#2}
} \providecommand{\href}[2]{#2}

\enddocument